\documentclass[12pt,reqno]{amsart}
\usepackage[margin=1.2in,marginparsep=0.1in,marginparwidth=1in]{geometry}
\usepackage{amssymb,amsmath,amsthm,amstext,amscd,latexsym,graphics,graphicx,caption, fullpage}
\usepackage[dvipsnames,svgnames,table]{xcolor}
\usepackage[plainpages=false,colorlinks=true, pagebackref]{hyperref}
\usepackage{tikz}
\usepackage{mathrsfs,enumerate,cleveref}
\usepackage{mathtools}
\usepackage{tikz-cd}
\hypersetup{citecolor=Sepia,linkcolor=blue, urlcolor=blue}
\usepackage{natbib}

\usepackage{bm}
\makeatletter
\def\@tocline#1#2#3#4#5#6#7{\relax
 \ifnum #1>\c@tocdepth 
 \else
 \par \addpenalty\@secpenalty\addvspace{#2}%
 \begingroup \hyphenpenalty\@M
 \@ifempty{#4}{%
  \@tempdima\csname r@tocindent\number#1\endcsname\relax
 }{%
  \@tempdima#4\relax
 }%
 \parindent\z@ \leftskip#3\relax \advance\leftskip\@tempdima\relax
 \rightskip\@pnumwidth plus4em \parfillskip-\@pnumwidth
 #5\leavevmode\hskip-\@tempdima
 \ifcase #1
 \or\or \hskip 2em \or \hskip 2em \else \hskip 3em \fi%
 #6\nobreak\relax
 \dotfill\hbox to\@pnumwidth{\@tocpagenum{#7}}\par
 \nobreak
 \endgroup
 \fi}
\makeatother
\newtheorem{intro-thm}{Theorem}[]
\theoremstyle{plain}
\newtheorem{thm}{Theorem}[section]
\newtheorem{theorem}[thm]{Theorem}
\newtheorem{question}[thm]{Question}
\newtheorem{lemma}[thm]{Lemma}

\theoremstyle{definition}

\newtheorem{rmks}[thm]{Remark}

\newcommand{\Det}{\mathrm{Det}}
\newcommand{\Tr}{{\rm Tr}}

\newcommand{\GL}{{\mathrm GL}}
\newcommand{\M}{\mathrm{M}}

\newcommand{\irr}{\mathrm{Irr}}

\newcommand{\sA}{{\mathcal A}}
\newcommand{\sB}{{\mathcal B}}

\newcommand{\sH}{{\mathcal H}}

\newcommand{\sN}{{\mathcal N}}
\newcommand{\sO}{{\mathcal O}}

\theoremstyle{remark}

\newcommand{\FF}{\mathbb{F}}

\usepackage{graphicx} 

\input{xy}

\title{Matrix Waring Problem in finite fields with coefficients from $\GL(2, \FF_p)$ }
 \author[S Datta]{Sutirtha Datta} \address{Department of Mathematics, IISER Pune, Dr Homi Bhabha Rd, Pashan, Pune, 411008, India} \email{sutirtha2702@gmail.com}
 \author[S Ghosh]{Sourav Ghosh} \address{Department of Mathematics, IISER Pune, Dr Homi Bhabha Rd, Pashan, Pune, 411008, India} \email{sourav.ghosh@students.iiserpune.ac.in}
 \author[A Maji]{Aniruddha Maji} \address{Department of Mathematics, IISER Pune, Dr Homi Bhabha Rd, Pashan, Pune, 411008, India} \email{maji.aniruddha@students.iiserpune.ac.in}
 \author[A Manik]{Antara Manik} \address{Department of Mathematics, IISER Pune, Dr Homi Bhabha Rd, Pashan, Pune, 411008, India} \email{manikantara66@gmail.com}

\begin{document}

 \thanks{Keywords: Waring problem, matrix over finite fields, non-abelian Fourier transform  }
 \subjclass[2020]{Primary 11P05, 15A24, 15B33 ; Secondary 20C15, 43A65}
 \begin{abstract}
   Let $k$ be a fixed positive integer and $p$ be an odd prime such that $\gcd(p-1,k)=1$. For sufficiently large primes $p$, we show that equation $AX^k+BY^k=C$ has solutions $X,Y \in \GL(2, \FF_p)$, where $A$ and $B$ are fixed matrices in $\GL(2, \FF_p)$, and $C$ is a fixed matrix in $\M(2,\FF_p).$ Our proof relies on Fourier analysis over $\GL(2, \FF_p).$
 \end{abstract}
	
 \maketitle
	\tableofcontents
\section{Introduction}
The classical Waring problem asks whether for every natural number $k$, there exists a positive integer $g(k)$ such that every natural number can be expressed as a sum of at most $g(k)$ many $k$-th powers. For example, $g(2)=4, g(3) =9, g(4)=19.$ This problem was first solved by Hilbert in $1909$ and later the formula for $g(k)$ was computed for all but finitely many $k.$ In a series of works, Hardy and Littlewood gave the definitive approach to the Waring problem, leveraging techniques from Fourier analysis. For a brief survey, refer to \cite{amm1971article}.

Modern versions consider a similar question over finite fields and objects with non-commutative structures. The matrix Waring problem asks whether every matrix over a commutative ring with unity $R$ can be expressed as a sum of two $k$-th powers of matrices over $R$. We are interested in the case when $R$ is the finite field $\FF_p$, where $p$ is an odd prime. \cite{kishore2022matrix} and \cite{kishore2025matrix} gave affirmative answers to the matrix Waring problem in finite fields, using Lang-Weil's results on the number of solutions to equations over finite fields. For recent developments, see \cite{kishore2026matrixwaringproblemiii}. 
 
We allow invertible coefficients and ask the following question:
\begin{question}\label{question}
    For a fixed positive integer $k$ and fixed invertible $2 \times 2$ matrices $A$ and $B$ over $\FF_p$, can a given $2 \times 2$ matrix $C$ with coefficients in $\FF_p$ be expressed as $$AX^k+BY^k=C?$$ 
\end{question}
where $X$ and $Y$ are invertible $2 \times 2$ matrices over $\FF_p$.
We set up some notations. 
Let $k\geq2$ be an integer, and consider the matrix algebra $M(2,\FF_p)$ and its unit group $\GL(2,\FF_p)$. Let us denote by $\sH$ the set of $k$-th powers in $\GL(2, \FF_p)$, where $d:=\gcd(k,p-1)=1$. Then Question \ref{question} asks for solution pairs $(X,Y) \in \sH \times \sH$ of the equation $X +BY=C$, where $B$ is a fixed matrix in $\GL(2, \FF_p)$ and $C$ is a fixed matrix in $\M(2, \FF_p).$ We give an affirmative answer for this question when $B \in \GL(2, \FF_p)$. Our proof relies on Fourier analysis over $\GL(2, \FF_p).$  

We recall the idea of the classical Gauss sum which blends additive and multiplicative characters. See \cite[Section 8.5]{apostol2013introduction} for a detailed treatment. Let $\psi$ be an additive character of $\FF_p$, and let $\chi$ be a multiplicative character of $\FF_p^{\times}.$ Extending $\chi$ to $\FF_p$ by setting $\chi(0)=0$, one defines $$G(\psi, \chi) = \sum_{x \in \FF_p} \psi(x) \chi(x)$$

as the Gauss sum. In the non-commutative setting, the analogous notion of matrix-valued Gauss sums for $\GL(n, \FF_p)$ was introduced by Lamprecht \cite{lamprecht1957struktur}. We fix $n=2$ and treat $\M(2, \FF_p)$ as an additive group so that the irreducible characters are given by $\psi_M(X):= \exp \left( \frac{2 \pi i }{p} \Tr(MX) \right).$ For each nonzero matrix $M \in\M(2, \FF_p)$, there exists a matrix $X$ such that $\Tr(MX) \neq 0$. Hence, the map $M \mapsto \psi_M$ is injective and yields an isomorphism of $\M(2,\FF_p)$ with its dual $\widehat{\M(2, \FF_p)}$. Then the non-abelian Gauss sum is defined as 

$$G(\rho ,A) : = \sum_{X \in \GL(2, \FF_p)} \rho(X)\psi_A(X)$$

where $\rho \in\irr(\GL(2, \FF_p))$ and $A \in \M(2, \FF_p)$. Ganguly and Rajan computed the trace values of the matrix Gauss sums for the irreducible representations of $\GL(2,\FF_p)$ in \cite{ganguly2023singular}. These will be explicitly stated wherever needed in the coming sections. Our Main result is the following:
\begin{theorem}\label{thm: main result}
   Let $k$ be a fixed positive integer, $p$ is an odd prime such that $\gcd(k, p-1) =1$, $B$ be a fixed matrix in $\GL(2,\FF_p)$ and $C$ be a fixed matrix in $\M(2,\FF_p)$. Then, for sufficiently large prime $p$, the equation $X^k+BY^k=C$ has solution $X$ and $Y$ in $\GL(2, \FF_p).$
\end{theorem}
When $B$ is scalar, this recovers \cite[Theorem A.6]{panjasainisingh25}, with the solutions actually in $\GL(2,\FF_p)$. 

\textbf{Layout of the Article:} In \Cref{sec: prelim}, we briefly recall the classification of conjugacy classes, inspect the subset of $k$-th power invertible matrices and determine its cardinality, and conclude with a brief discussion on the character values of irreducible representations of $\GL(2,\FF_p)$ and some useful lemmas for later use. In \Cref{sec: main idea of proof}, we record the key idea and define the function $S_k(M)$ in terms of the Fourier coefficients, which we compute in \Cref{sec: computing Fourier coeff} for all four kinds of irreducible representations. Finally in \Cref{sec: proof of the main theorem}, we estimate $S_k(M)$ for the zero matrix, non-singular and rank one matrices in subsections \ref{subsec: zero matrix}, \ref{subsec: nonsingular} and \ref{subsec:rank one} before concluding the proof of the main theorem.

\subsection*{Acknowledgement} The authors thank Anupam Singh and Steven Spallone for several helpful comments.  

\section{Preliminaries on representation theory and structure of $\GL(2,\FF_p)$}\label{sec: prelim}
\subsection{Conjugacy classes and $k$-th power matrices in $\GL(2, \FF_p)$}
We compute the size of $\sH$ using the conjugacy classes of $\GL(2, \FF_p)$. We briefly recall the four types of conjugacy classes in $\GL(2, \FF_p)$. First, there are the split semisimple elements, which are diagonalizable matrices with distinct eigenvalues. The size of this conjugacy class is $p(p+1).$ Such a matrix will be a $k$-th power iff both the eigenvalues are $k$-th powers in $\FF_p^{\times}$. The number of $k$-th powers in $\FF_p^{\times}$ is $p-1$, as $\gcd(k, p-1)=1.$ Hence, the total number of $k$-th power matrices of first type is $\binom{p-1}{2} (p^2+p).$
 
Secondly, there are elliptic elements whose characteristic polynomial is irreducible over $\FF_p$, possessing conjugate eigenvalues in $\FF_{p^2}^{\times }  \setminus \FF_p^{\times}$. Fixing a non-square $\varepsilon\in \FF_p^{\times}$, every such class is represented by $\begin{pmatrix}
a & b\varepsilon\\ b & a
\end{pmatrix}$
for some $a\in \FF_p$ and $b\in \FF_p^{\times}$. The corresponding eigenvalues are $\lambda:=a+b\sqrt{\varepsilon}$ and its conjugate $\lambda^p=a-b\sqrt{\varepsilon}$ in $\FF_{p^2}^{\times}$. Thus, each conjugacy class of this type contains exactly $p(p-1)$ elements. Such a matrix will be a $k$-th power iff $\lambda$ is a $k$-th power in $\FF_{p^2}^{\times} \setminus \FF_p^\times.$ The number of $k$-th powers in $\FF_{p^2}^{\times} \setminus \FF_p^\times$ is $\frac{p^2-1}{d'} - (p-1)$, where $d' = \gcd(k, p^2-1).$ Thus, $\lambda$ has $\frac{p^2-1}{d'} - (p-1)$ many choices as $\lambda \notin \FF_p$. Since $\lambda \neq \lambda^p$, the number of pairs $\{\lambda, \lambda^p\}$ is $\frac{1}{2}\left(\frac{p^2-1}{d'} - (p-1)\right)$. So the total number of $k$-th power matrices of the second type is $\frac{1}{2}\left(\frac{p^2-1}{d'} - (p-1)\right) (p^2-p).$
 
Thirdly, there are the scalar matrices $\begin{pmatrix}
    a & 0 \\ 0 & a
\end{pmatrix}$. There are $p-1$ such matrices and all the matrices are $k$-th powers.
Finally, we have the non-semisimple elements of the form $\begin{pmatrix}
    a & 1 \\ 0 & a
\end{pmatrix}$, whose conjugacy class has size $p^2-1.$ The total number of such matrices which are $k$-th power is $(p^2-1)(p-1).$ We add the count to get  
\begin{align*}
    |\sH| &= \binom{p-1}{2} (p^2+p) + \frac{1}{2}\left(\frac{p^2-1}{d'} - (p-1)\right) (p^2-p)+ p-1+ (p^2-1)(p-1)\\
    &=p(p-1)(p^2-1) \left(\frac{1}{2} + \frac{1}{2d'}\right) .
\end{align*} 
\subsection{Character Values of Irreducible Representations of $\GL(2,\FF_p)$}
The four types of irreducible characters and their values on different conjugacy classes is summarized below.
\[
\begin{tabular}{ |c|c|c|c|c| }
 \hline
$A$ & $\begin{pmatrix}
    a & 0\\ 0&a
\end{pmatrix}$ & $\begin{pmatrix}
    a & 0\\ 0&d
\end{pmatrix}, a\neq d$ & $\begin{pmatrix}
    a &1 \\ 0&a
\end{pmatrix}$ & $\begin{pmatrix}
    a &b \varepsilon \\ b&a
\end{pmatrix}, b\neq0$\\\hline\hline
   $\mu \circ \det$ & $\mu(a)^2 $&$\mu(a)\mu(d)$ &$\mu(a)^2$ & $\mu(a^2-b^2\varepsilon)$ \\
   \hline
    $\mathrm{St}\otimes \chi$ & $p\chi(a)^2$ & $\chi(a)\chi(d)$ &$0$&$-\chi(a^2-b^2\varepsilon)$\\
    \hline
  $I_{\chi_1,\chi_2}$ & $(p+1)\chi_1(a)\chi_2(a)$ & $\chi_1(a)\chi_2(d)+\chi_1(d)\chi_2(a)$ &$\chi_1(a)\chi_2(a)$& $0$ \\
  \hline
  $\Theta_{\phi}$ & $(p-1)\phi(a)$ & $0$ &$-\phi(a)$& $-(\phi(\lambda)+\phi(\lambda^p))$\\\hline
\end{tabular}
\]
\subsection{Some Important Lemmas}
Let $p$ be an odd prime. Denoting the set of nonzero squares as $S=\{t^2:t\in \FF_p^{\times}\}$, the Legendre symbol of $a\in \FF_p$ is given by 
$$\left(\frac{a}{p}\right)=\begin{cases}
    0 & \text{~when~} a=0\\
    1 & \text{~when~} a\in S\\
    -1 & \text{~when~} a\in \FF_p^{\times}\setminus S
\end{cases}
$$
For coefficients $a,b,c\in \FF_p$, let $\Delta=b^2-4ac$ denote the discriminant of the quadratic $ax^2+bx+c$.
\begin{lemma}\label{lem: legendre lemma} If $a\neq 0$ and $b,c, \Delta$ are as described above, then  
$$\sum_{x\in \FF_p}\left(\frac{ax^2+bx+c}{p}\right)=\begin{cases}
    (p-1)\left(\frac{a}{p}\right) & \text{~when~}\Delta=0\\
    -\left(\frac{a}{p}\right) & \text{~when~}\Delta\neq 0
\end{cases}
$$

\end{lemma}
\begin{proof}
Completing the square yields $4a(ax^2+bx+c)=(2ax+b)^2-\Delta$. By the multiplicative property of the Legendre symbol, we factor out $\left(\frac{a}{p}\right)$ using $\left(\frac{4a^2}{p}\right)=1$, which under the linear bijection $u=2ax+b$ on $\FF_p$ reduces the sum to 
$$\sum_{x\in \FF_p}\left(\frac{ax^2+bx+c}{p}\right)=\sum_{x\in \FF_p}\left(\frac{4a^2}{p}\right)\left(\frac{ax^2+bx+c}{p}\right)=\left(\frac{a}{p}\right)\sum_{u\in \FF_p}\left(\frac{u^2-\Delta}{p}\right).
$$
Applying the monic quadratic evaluation $\sum_{u\in \FF_p}\left(\frac{u^2-\Delta}{p}\right)=p-1$ for $\Delta=0$ and $-1$ for $\Delta\neq 0$ \cite[Exercise~ $6$]{ireland1990rosen} completes the proof.
\end{proof}
 \begin{lemma}\label{lem: elliptic M goes to scalar MB iff B is nonsplit ss}
    For an elliptic $M$, $MB$ is scalar only if $B$ is elliptic.
\end{lemma}
\begin{proof}
    Assume that $MB=cI$, for some $c \in \FF_p^{\times}$. Then $B= c M^{-1}$, which is elliptic.
\end{proof}
 
For a given $B \in \GL(2, \FF_p),$ define $\alpha(B)$ to be the number of elliptic matrices $M$ such that $MB$ is elliptic. $MB$ is elliptic if and only if the discriminant $\Delta_{MB}$ is non-square in $\FF_p^{\times}.$ Let $c\in \FF_p^{\times}$. Then  
$$\Delta_{M(cB)}=\Tr(cMB)^2 - 4 \Det(cMB) = c^2 \left(\Tr(MB)^2 -4\Det(MB)\right)=c^2 \Delta_{MB}$$
Therefore, $\Delta_{M(cB)}$ is non-square in $\FF_p^{\times}$ if and only if $\Delta_{MB}$ is non-square. Consequently, $\alpha$ is invariant under constants: $\alpha(B) =\alpha(cB)$ for all $c \in \FF_p^{\times}.$ Moreover, $\alpha$ is also conjugacy invariant. Indeed, if $B':= g^{-1}Bg$, then $MB'= g^{-1}gMg^{-1}Bg$. Denote $M'=gMg^{-1}.$ Thus, $MB'$ and $M'B$ are conjugate, which implies that $MB'$ is elliptic if and only if $M'B$ is elliptic. Now, $\alpha(B')$ is the number of elliptic $M$ such that $MB'$ is elliptic. In other words, it is the number of elliptic $M'$ such that $M'B$ is elliptic, since the conjugation on $M$ is a bijection of elliptic elements. 

For $D \in \FF_p^{\times},$ we define $\beta(D):=|\sB| $ where $\sB:=\{ t \in \FF_p \mid t^2 -4D \text{~is non-square in~} \FF_p \}$. Denote 
$\sB_s:=\{ t \in \FF_p \mid 0 \neq t^2 -4D \text{~is square in~} \FF_p \}$ and $\sB_0:=\{ t \in \FF_p \mid t^2 -4D=0 \}$. We have 
 \begin{align*}
    &\sum_{t\in \FF_p}\left(\frac{t^2-4D}{p}\right) = \sum_{t \in \sB}\left(\frac{t^2-4D}{p}\right)+\sum_{t \in \sB_s}\left(\frac{t^2-4D}{p}\right) +\sum_{t \in \sB_0}\left(\frac{t^2-4D}{p}\right)=|\sB_s|-|\sB|
\end{align*}
Since $|\sB|+ |\sB_s| +|\sB_0| = p,$ we have $2|\sB|= p-|\sB_0| - \sum_{t\in \FF_p}\left(\frac{t^2-4D}{p}\right)= p-\left(\frac{D}{p}\right).$
Thus 

\begin{equation}\label{eq: beta D values}
    \beta(D)=\frac{1}{2}\left(p-\left(\frac{D}{p}\right)\right)=\begin{cases}
    \frac{p+1}{2} & \text{~when~} \left(\frac{D}{p}\right)=-1\\
    \frac{p-1}{2} & \text{~when~} \left(\frac{D}{p}\right)=1 
    \end{cases}
\end{equation}

\begin{lemma}\label{lem: nonsplit ss M goes to nonsplit ss MB}
   Suppose $B$ is a non-scalar matrix in $\GL(2, \FF_p)$. Then $$\alpha  (B) \leq \frac{p^4-1}{4}.$$  
\end{lemma}
\begin{proof} 
    First, assume that $B$ is a non-semisimple element. By the constant invariance and conjugacy invariance of $\alpha$, we can take $B:=\begin{pmatrix}
        1 & 1\\ 0&1
    \end{pmatrix}$. Let $M =\begin{pmatrix}
         a & b \\ c&d
    \end{pmatrix},$ so that $MB= \begin{pmatrix}
        a& a+b \\ c & c+d
    \end{pmatrix}.$ The conjugacy class of elliptic elements is determined uniquely by its trace and determinant. Suppose $\Tr(M)=t_1, \Det(M)=D, \Tr(MB)= t_2.$ It follows that $\Det(MB)=D.$ Fixing $t_1, t_2, D$ and $a$ determines $M$ entirely. Hence, the number of such $M$ is $p.$ There are $\beta(D)$ many choices for $t_1$ such that $t_1^2-4D$ is non-square. Since $c=t_2-t_1\neq 0$, there are $\beta(D) -1$ many choices for $t_2$ such that $t_2^2 - 4D $ is non-square. Hence, for a fixed $D$, there are $p\beta(D)(\beta(D)-1)$ many choices. Hence, 
    \begin{align*}
        \alpha(B)&=\sum_{D \text{~square}} p\beta(D)(\beta(D)-1) +\sum_{D \text{~non-square}} p\beta(D)(\beta(D)-1)\\
        &= \sum_{D \text{~square}} p\frac{p-1}{2}\frac{p-3}{2} +\sum_{D \text{~non-square}} p\frac{p+1}{2}\frac{p-1}{2} \\
        &=\frac{p(p-1)^3}{4}
    \end{align*}
    Next, assume that $B=\begin{pmatrix}
        \lambda_1 & 0\\ 0&\lambda_2
    \end{pmatrix}$ and $M=\begin{pmatrix}
        a &b \\c &d 
    \end{pmatrix}$. Then, $MB= \begin{pmatrix}
        a \lambda_1 &  b\lambda_2\\ c \lambda_1 & d\lambda_2
    \end{pmatrix}.$ Suppose $\Tr(M)=t_1, \Det(M)=D_1, \Tr(MB)= t_2$ and $\Det(MB)=D_2=D_1\lambda_1\lambda_2.$ Fixing $t_1, t_2, D_1, B,b$ completely determines $M$. Hence, the number of such $M$ is $p-1.$ There are $\beta(D_1)$ many choices for $t_1$ such that $t_1^2-4D_1$ is non-square. There are $\beta(D_2)$ many choices for $t_2$ such that $t_2^2 - 4D_2 $ is non-square. Hence 
    \begin{align*}
        \alpha(B) &=  \sum_{D_1, D_2 \in \FF_p^{\times}} (p-1)\beta(D_1)\beta(D_2) 
    \end{align*}
    When $\lambda_1\lambda_2$ is a square, $\left(\frac{D_1}{p}\right)=\left(\frac{D_2}{p}\right) $. Thus,
    \begin{align*}
        \alpha(B) &= \sum_{D_1 \text{~square}}   (p-1)\beta(D_1)\beta(D_2)+\sum_{ D_1 \text{~non-square}} (p-1)\beta(D_1)\beta(D_2)= \frac{(p-1)^2(p^2+1)}{4}
    \end{align*}
    When $\lambda_1\lambda_2$ is non-square, $\left(\frac{D_1}{p}\right)=-\left(\frac{D_2}{p}\right) $. Thus,
    \begin{align*}
        \alpha(B) &= \sum_{D_1 \text{~square}}   (p-1)\beta(D_1)\beta(D_2)+\sum_{ D_1 \text{~non-square}} (p-1)\beta(D_1)\beta(D_2)=\frac{(p-1)^2(p^2-1)}{4}
    \end{align*}
    Finally, assume that $B$ is elliptic, in terms of the companion matrix, say $B= \begin{pmatrix}
        0 & -D_B\\ 1 &t_B
    \end{pmatrix}$, $\Tr(B)=t_B$ and $\Det(B)=D_B$ and $M=\begin{pmatrix}
        a &b \\ c&d
    \end{pmatrix}$,  Then $MB= \begin{pmatrix}
        b & bt_B -aD_B\\ d  & dt_B-cD_B
    \end{pmatrix}$. Suppose  $\Tr(M)=t_1, \Det(M)=D_1, \Tr(MB)= t_2,$ and $\Det(MB)=D_2.$ Note that $D_2=D_BD_1, t_1=a+d$ and $t_2= b+dt_B-cD_B$. From the last two equations we have $b = cD_B + t_2 + t_B(a-t_1)$. Since $bc-ad=D_1$, we also have 
    $$c^2D_B + c(t_2+t_B(a-t_1))=at_1-a^2-D_1$$ 
    Let $f(c):=c^2D_B + c(t_2+t_B(a-t_1))-at_1+a^2+D_1$ be a polynomial in $c.$ For a fixed $a \in \FF_p$, the discriminant $\Delta_f=(cD_B+t_B(a-t_1))^2-4D_B(D_1+a^2-at_1)$ is a polynomial in $a$. If $\Delta_f=0$, then there is exactly one solution for $c.$ When $\Delta_f$ is non-square, then there is no such $c$. For a nonzero square $\Delta_f$, there are two solutions to $f(c)=0.$ Fixing $D_1, B, t_1, t_2, a$ and $c$ completely determines $M$. Consequently, the total number of elliptic $M$ such that $MB$ is elliptic is $$\sum_{a\in \FF_p} \left(1+ \left(\frac{\Delta_f}{p}\right)\right)= p+\sum_{a\in \FF_p}   \left(\frac{\Delta_f}{p}\right) = p- \left(\frac{\Delta_B}{p}\right) ,$$
    where the last equality follows by \Cref{lem: legendre lemma}, since $\Delta_B$ is non-square. Hence, we conclude that there are at most $p+1$ many such $M$. Therefore,
        \begin{align*}
        \alpha(B) &=  \sum_{D_1\in \FF_p^{\times}} (p+1)\beta(D_1)\beta(D_2)=\begin{cases}
             \frac{p^4-1}{4} & \text{~when~}D_B \text{~is square}\\
            \frac{(p^2-1)^2}{4} & \text{~when~} D_B \text{~is non-square}
        \end{cases}  
    \end{align*}
\end{proof}
For a given $B$ define $\gamma(B)$ to be the number of elliptic matrices $M$ such that $MB$ is non-semisimple.
\begin{lemma}\label{lem: nonsplit ss M goes to non-semisimple MB}
   Suppose $B$ is a non-scalar matrix in $\GL(2, \FF_p)$. Then $$\gamma(B)\leq \frac{(p+1)(p-1)^2}{2} .$$  
\end{lemma}
\begin{proof}
    First, we assume that $B = \begin{pmatrix}
        \lambda_1 & 0\\ 0 & \lambda_2
    \end{pmatrix}$ and $M=\begin{pmatrix}
        a & b\\ c & d
    \end{pmatrix}$. Denote $\Det(M)$ by $D$ and $\Tr(M)$ to be $t.$ Then $MB = \begin{pmatrix}
        a\lambda_1 & b \lambda_2\\
        c\lambda_1 & d \lambda_2
    \end{pmatrix}.$ Since $MB$ is non-semisimple, it is conjugate to $\begin{pmatrix}
        x & 1\\ 0 &x
    \end{pmatrix}$ for some $x \in \FF_p^{\times}.$ Thus, $\Tr(MB)=2x=a\lambda_1+d\lambda_2$ and $\Det(MB)=x^2=\lambda_1\lambda_2D$. Fixing $t, D, x,B$ and $b$ determines the matrix $M$, giving $(p-1)$-many choices of $M$ and two choices for $x$. There are $\beta(D)$ many elements such that $t^2-4D$ is non-square. Therefore,
    \begin{align*}
        \gamma(B) =\begin{cases}
            \sum\limits_{D \text{~square}} 2(p-1)\beta(D) &=\frac{(p-1)^3}{2} \text{~~when~} \left(\frac{\lambda_1\lambda_2}{p}\right)=1 \\
            \sum\limits_{D \text{~non-square}} 2(p-1)\beta(D)&=\frac{(p+1)(p-1)^2}{2} \text{~~when~} \left(\frac{\lambda_1\lambda_2}{p}\right)=-1
        \end{cases} 
    \end{align*}
    Next, we assume that $B=\begin{pmatrix}
        y & 1\\ 0 &y
    \end{pmatrix}$ and $M=\begin{pmatrix}
        a & b\\ c & d
    \end{pmatrix}.$ Denote $\Det(M)$ by $D$ and $\Tr(M)$ to be $t.$  Then $MB = \begin{pmatrix}
        ay & a+by\\
        cy & c+dy
    \end{pmatrix}.$ Again, assume that $MB$ is conjugate to $\begin{pmatrix}
        x & 1\\ 0 &x
    \end{pmatrix}$ for some $x \in \FF_p^{\times}.$ Thus, $\Tr(MB)=2x=ay+dy+c$ and $\Det(MB)=x^2=y^2D$, implying that $D$ is square. Fixing $t, D, x,B$ and $b$ determines the matrix $M$, giving $(p-1)$-many choices of $M$ and two choices for $x$. Therefore,
    \begin{align*}
        \gamma(B)= \sum_{D  \text{~square}}   2(p-1)\beta(D) = \frac{(p-1)^3}{2} 
    \end{align*}
     Finally, assume that $B$ is elliptic, in terms of the companion matrix, say $B= \begin{pmatrix}
        0 & -D_B\\ 1 &t_B
    \end{pmatrix}$, $\Tr(B)=t_B$ and $\Det(B)=D_B$ and $M=\begin{pmatrix}
        a &b \\ c&d
    \end{pmatrix}$,  Then $MB= \begin{pmatrix}
        b & bt_B -aD_B\\ d  & dt_B-cD_B
    \end{pmatrix}$. Suppose  $\Tr(M)=t_1, \Det(M)=D_1, \Tr(MB)= t_2,$ and $\Det(MB)=D_2.$ Note that $D_2=D_BD_1, t_1=a+d$ and $t_2= b+dt_B-cD_B$. Suppose $MB$ is conjugate to $\begin{pmatrix}
        x & 1\\ 0 &x
    \end{pmatrix}$ for some $x \in \FF_p^{\times}.$ Similar to the last case of the proof of \Cref{lem: nonsplit ss M goes to nonsplit ss MB}, we have a quadratic polynomial in $c$ $$c(2x + (a-t_1)t_B +c D_B) =a(t_1-a)-D_1$$ Therefore, the number of elliptic matrices $M$ such that $MB$ is nonsemisimple is at most $p+1.$ Now $B=xM^{-1}$ as of \Cref{lem: elliptic M goes to scalar MB iff B is nonsplit ss} leads to $MB=\begin{pmatrix}
        x &0\\ 0 &x
    \end{pmatrix}$, which shares the same trace and determinant as $\begin{pmatrix}
        x &1\\ 0 &x
    \end{pmatrix}$. Therefore, we exclude the corresponding choice of $t_1$ from the count of $\beta(D_1).$ Also $x^2=D_2$ gives two choices for $x$ leading to
    \begin{align*}
        \gamma(B) =\begin{cases}
            \sum\limits_{D_1 \text{~square}} 2(p+1)(\beta(D_1)-1)&=\frac{(p-3)(p^2-1)}{2} \text{~~when~} \left(\frac{D_B}{p}\right)=1 \\
            \sum\limits_{D_1 \text{~non-square}} 2(p+1)(\beta(D_1)-1)&=\frac{(p+1)(p-1)^2}{2} \text{~~when~} \left(\frac{D_B}{p}\right)=-1
        \end{cases} 
    \end{align*}
\end{proof}
Let us define $A_a:= \begin{pmatrix}
    a & 0 \\
    0 &0
\end{pmatrix}$, $a\neq 0$ and $N:= \begin{pmatrix}
    0 &1\\
    0 &0
\end{pmatrix}$. Then every rank one matrix in $\M(2, \FF_p)$ is conjugate to either $A_a$ or $N.$
\begin{lemma}
    There are exactly $(p^2-1)$ many elements conjugate to $N$ and $p(p^2-1)$ many elements conjugate to $A_a$ in $\M(2, \FF_p)$.
\end{lemma}
\begin{proof}
    This follows by orbit-stabilizer theorem.
\end{proof}

\begin{lemma}\label{lem: A_a goes to N and vice versa}
    Suppose $B$ is a non-scalar matrix in $\GL(2, \FF_p)$. Then for at most $p^2-1$ many $M$ conjugate to $N$, $MB$ is conjugate to $A_a$ and vice versa.
\end{lemma}
\begin{proof}
    Since $M \sim N$, $M$ must have both the trace and determinant $0$ and hence, it will be of the form $\begin{pmatrix}
        x & y \\ z & -x
    \end{pmatrix}$ such that $x^2 +yz=0.$
    
    \textbf{Case-$1$} When $B= \begin{pmatrix}
        b & 1\\
        0 & b
    \end{pmatrix}$, we have $\Tr(MB)=z$. For $MB$ to be conjugate to $A_a$, we must have $z=a.$ Then $$M =\begin{pmatrix}
        x & -x^2a^{-1}\\ a & -x
    \end{pmatrix}$$ For a fixed $a\in\FF_p^{\times}$, we have $p$ many choices of $M$. Since $a\neq 0,$ total number of $M$ is $p(p-1).$
    
    \textbf{Case-$2$} When $B= \begin{pmatrix}
        \lambda_1 & 0\\
        0 & \lambda_2
    \end{pmatrix}$, we have $\Tr(MB)=x(\lambda_1 - \lambda_2)$, which must be equal to $a$, if $MB$ has to be conjugate to $A_a.$ Then 
    $$M = \begin{pmatrix}
        a(\lambda_1- \lambda_2)^{-1} & y\\ -a^2y^{-1}(\lambda_2- \lambda_1)^{-2} & a(\lambda_2- \lambda_1)^{-1}
    \end{pmatrix}$$
    For a fixed $a\in\FF_p^{\times}$, the number of such matrices $M$ is $p-1$, since $a \neq 0$ and $\lambda_1 \neq \lambda_2$ implies $y$ to be nonzero. This makes the total count in this case to be $(p-1)^2.$
    
    \textbf{Case-$3$} When $B$ is elliptic, say $\begin{pmatrix}
        0 & -D_B\\ 1 &t_B
    \end{pmatrix}$, we have $\Tr(MB)=y-D_Bz-xt_B=a$, as $MB$ is conjugate to $A_a.$ Hence, $y=a+D_Bz+xt_B$. By the similar argument as the elliptic $B$ cases in the previous lemmas, there are $p+1$ choices for $M.$ Since $a$ has $p-1$ choices, the total count is $p^2-1.$
\end{proof}
\begin{rmks}\label{A_a goes to A_a and N goes to N}
    Suppose $B$ is non-scalar. Then for at most $p(p^2-1)-(p-1)^2$ many $M$ conjugate to $A_a$, $MB$ is conjugate to $A_a$. Similarly,  at most $(p^2-1)-(p-1)^2$ many $M$ conjugate to $N$, $MB$ is conjugate to $N$. We will use these in the computation of the rank one case.
\end{rmks}
\section{Main Idea of the Proof}\label{sec: main idea of proof}
Let us denote the number of solutions $(X,Y) \in \sH \times \sH$ to $X+BY=C$ as $\sN(C).$ Observe that $$\frac{1}{p^4}\sum_{M \in \M(2, \FF_p)}\psi_M(A)= \begin{cases}
    1 \quad \text{when~} A =0\\
    0 \quad \text{when~} A \neq 0.
\end{cases}$$ by character orthogonality relations. Then we have  
\begin{align*}
    \sN (C) &= \sum_{X,Y \in H}\left( \frac{1}{p^4} \sum_{M \in \M(2, \FF_p)}\psi_M(X+BY-C)\right)\\ &=  \frac{1}{p^4} \sum_{M \in \M(2, \FF_p)} \sum_{X,Y \in H}\psi_M(X)  \psi_{MB} (Y)  \psi_M(-C)
\end{align*}
where we used the linearity and permutation invariance of matrix trace. Define $$S_k(M):= \sum\limits_{X \in \sH} \psi_M(X),$$ which allows us to rewrite $\sN(C)$ as 
\begin{equation}\label{eq: N(C) formula}
    \sN(C)= \frac{1}{p^4} \sum_{M \in \M(2, \FF_p)} \psi_M(-C)S_k(M) S_k(MB)
\end{equation}
Using the indicator function $1_{\sH}$ we extend $S_k(M)$ to $\GL(2, \FF_p)$ : $$S_k(M) = \sum_{X \in \GL(2, \FF_p)} 1_{\sH}(X)\psi_M(X)$$
Since $1_{\sH}$ is a class function, it can be written as a linear combination of the irreducible characters $$1_{\sH} (X)= \sum_{\rho \in \irr(\GL(2, \FF_p))}  c_{\rho} \chi_{\rho}(X)$$
Thus $S_k(M)$ can be written as 
\begin{equation}
\label{Eq:S_k(M)} S_k(M) = \sum_{X \in \GL(2, \FF_p)}\sum_{\rho \in \irr(\GL(2, \FF_p))}  c_{\rho} \chi_{\rho}(X)\psi_M(X)= \sum_{\rho\in \irr(\GL(2, \FF_p))}c_{\rho}\Tr(G(\rho,M)) 
\end{equation}
\section{Fourier Coefficients Computations}\label{sec: computing Fourier coeff}
In this section, we compute the Fourier Coefficients $c_{\rho}$ for each of the four types of irreducible representations of $\GL(2, \FF_p).$
\subsection{$1$-Dimensional Representations} For the trivial representation $1_G$ we have $$c_{1_G}= \frac{1}{|G|} \sum\limits_{M \in \sH} 1_G(M)= \frac{|\sH|}{|G|}=\frac{1}{2}+\frac{1}{2d'}.$$
For the other one dimensional representations $\rho$, we compute $c_{\pi}$ by considering each of the cases when $M$ is  scalar, non-semisimple, split semisimple and elliptic. $H'$ will denote the subgroup of $k$-th powers in $\FF_{p^2}^{\times}$ and we will use the fact that $\sum_{\alpha \in H'}\overline{\chi'(\alpha)}=0$ for any nontrivial character $\chi'$ of $\FF_{p^2}^{\times}.$ the cardinality of $H'$ will be $\frac{p^2-1}{d'}$, where $d' = \gcd ( p^2-1,k).$
\begin{align*}
   c_{\rho} &= \frac{1}{|G|}(\underbrace{\sum_{a \in H} \overline{\chi(a^2)}}_{\text{scalar}} + \underbrace{\sum_{a \in H} (p^2-1) \overline{\chi(a^2)}}_{\text{non-semisimple}}+ \underbrace{\frac{p(p+1)}{2}\sum_{a\neq b \in H} \overline{\chi(a) \chi(b)}}_{\text{split semisimple}}\\&+ \underbrace{\frac{p(p-1)}{2} \sum_{\alpha, \alpha^p \in H'-H}\overline{\chi(\alpha^{p+1})}}_{\text{elliptic}})\\ 
   &= \frac{1}{|G|}\left(\sum_{a \in H} \overline{\chi(a^2)}+ \sum_{a \in H} (p^2-1) \overline{\chi(a^2)}+ \frac{p(p+1)}{2}\left(\sum_{a, b \in H} \overline{\chi(a) \chi(b)} -\sum_{a \in H} \overline{\chi(a^2)}\right)\right)\\ 
   &+ \frac{p(p-1)}{2|G|}\left(\sum_{\alpha, \alpha^p \in H'}\overline{\chi(\alpha^{p+1})} -\sum_{\alpha\in H}\overline{\chi(\alpha^2)}\right)  \\
   &= \frac{1}{|G|}\left(\sum_{a \in H} \overline{\chi(a^2)}+ \sum_{a \in H} (p^2-1) \overline{\chi(a^2)}+ \frac{p(p+1)}{2}\left(\sum_{a\in H} \overline{\chi(a)}\sum_{b\in H} \overline{\chi(b)} -\sum_{a \in H} \overline{\chi(a^2)}\right)\right)\\
   &+ \frac{p(p-1)}{2|G|}\left(\sum_{\alpha  \in H'}\overline{\chi(N_{\FF_{p^2}/ \FF_p}\alpha)} -\sum_{\alpha\in H}\overline{\chi(\alpha^2)}\right) \\
   & = \frac{1}{|G|}\left(\sum_{a \in H} \overline{\chi(a^2)}+ \sum_{a \in H} (p^2-1) \overline{\chi(a^2)}+ \frac{p(p+1)}{2} \left( -\sum_{a \in H} \overline{\chi(a^2)}\right)\right)\\
   &+ \frac{p(p-1)}{2|G|}\left(\sum_{\alpha \in H'}\overline{\chi'(\alpha)} -\sum_{a\in H}\overline{\chi(a^2)}\right) \\
   &=\frac{1}{|G|}\left(1+(p^2-1)- \frac{p(p+1)}{2}-\frac{p(p-1)}{2}\right)\sum_{a \in H} \overline{\chi(a^2)} = 0 
\end{align*}
 Note that $\chi':= \chi \circ N_{\FF_{p^2}/ \FF_p}$ is nontrivial since $\chi$ is nontrivial.
\subsection{Steinberg Representations}
 Suppose $\rho = \mathrm{St}\otimes \chi$, where $\chi$ is a nontrivial character of $\FF_p^\times$. Then we have 
 \begin{align*}
     & |G|c_{\rho} =  \underbrace{p \sum_{a \in H} \overline{\chi(a^2)}}_{\text{scalar}}+\underbrace{0}_{\text{non-semisimple}}+ \underbrace{\frac{(p^2+p)}{2} \sum_{a \neq d  \in H} \overline{\chi(a) \chi(d)}}_{\text{split semisimple}} + \underbrace{\frac{(p^2-p)}{2} \sum_{\{\alpha, \alpha^p\} \subseteq H'-H}  \overline{-\chi(\alpha^{p+1})}}_{\text{elliptic}} \\
     &=   p \sum_{a \in H} \overline{\chi(a^2)} +  \frac{(p^2+p)}{2}\left( \sum_{a  \in H} \overline{\chi(a)} \sum_{d \in H} \overline{\chi(d)}  - \sum_{a  \in H} \overline{\chi(a^2)} \right)\\ &+ \frac{(p^2-p)}{2} \left( \sum_{\alpha \in H'}  \overline{-\chi(\alpha^{p+1})}- \sum_{\alpha \in H}  \overline{-\chi(\alpha^2)}\right) \\
     &=0
 \end{align*}
 When $\chi$ is trivial, we have
\begin{align*}
      c_{\mathrm{St}} &=\frac{1}{|G|} \left( \frac{p^2+p}{2}\cdot |H|^2-\frac{p^2-p}{2}\cdot |H'|\right)= \frac{1}{|G|} \left(\frac{p^2+p}{2}\cdot (p-1)^2-\frac{p^2-p}{2}\cdot\frac{p^2-1}{d'}\right)\\&= \frac{1}{2}-\frac{1}{2d'}
\end{align*}
 \subsection{Principal Series} Suppose $\rho= I_{\chi_1, \chi_2}$, where $\chi_1$ or $\chi_2$ is nontrivial characters of $\FF_p^\times$. Then 
 \begin{align*}
      |G|c_{\rho} &=  \underbrace{(p+1) \sum_{a \in H} \overline{\chi_1(a)\chi_2(a)}}_{\text{scalar}} + \underbrace{(p^2-1) \sum_{a\in H} \overline{\chi_1(a)\chi_2(a)}}_{\text{non-semisimple}}\\&+ \underbrace{p(p+1)\sum_{\{a,b\} \subseteq H} \overline{(\chi_1(a)\chi_2(b) + \chi_1(b) \chi_2(a)) } }_{\text{split semisimple}} + \underbrace{0}_{\text{elliptic}}\\
     &= (p^2+p) \sum_{a\in H} \overline{\chi_1(a)\chi_2(a)} +  p(p+1)\sum_{a \neq b \in H} \overline{\chi_1(a)\chi_2(b)}   \\
     &=(p^2+p) \sum_{a\in H} \overline{\chi_1(a)\chi_2(a)} +  p(p+1)\left(\sum_{a \in H} \overline{\chi_1(a)} \sum_{b \in H} \overline{\chi_2(b)} - \sum_{a \in H} \overline{\chi_1(a)\chi_2(a)}   \right)\\ &=0\\
 \end{align*}
 \subsection{Cuspidal Representations} Suppose $\rho= \Theta_\phi$, where $\phi$ is a nontrivial character of $F_{p^2}^\times$ such that $\phi ^p\neq \phi$.
\begin{align*}
    |G| c_{\rho}&= \underbrace{(p-1) \sum_{a\in H}\overline{\phi(a)}}_{\text{scalar}} + \underbrace{(p^2-1) \sum_{a\in H} \overline{-\phi(a)} }_{\text{non-semisimple}}+ \underbrace{0}_{\text{split semisimple}} \\&+ \underbrace{(p^2-p) \sum_{\{\alpha, \alpha^p\} \subseteq H'-H}\overline{-(\phi(\alpha) + \phi(\alpha^p))} }_{\text{elliptic}}\\
    & = (p-1) \sum_{a\in H}\overline{\phi(a)} +  (p^2-1) \sum_{a\in H} \overline{-\phi(a)}  +    (p^2-p) \sum_{\alpha \in H'-H}\overline{-(\phi(\alpha))}\\
    & = (p-1) \sum_{a\in H}\overline{\phi(a)} +  (p^2-1) \sum_{a\in H} \overline{-\phi(a)}  +    (p^2-p) \left(\sum_{\alpha \in H'}\overline{-(\phi(\alpha))}- \sum_{\alpha \in H}\overline{-(\phi(\alpha))} \right)\\
    &= (p^2-p) \sum_{\alpha \in H'}\overline{-(\phi(\alpha))} + \left(p-1-p^2+1+p^2-p \right)\sum_{a\in H}\overline{\phi(a)}   \\
    & = (p^2-p) \sum_{\alpha \in H'}\overline{-(\phi(\alpha))}
\end{align*}
The last term is zero when $\phi^{d'}\neq 1$, where $d' := \gcd(p^2-1,k),$ as $\phi|_{H'}$ is nontrivial. For the remaining $\frac{d'-1}{2}$ many characters $\phi$ such that $\phi^{d'}=1$, we have $\phi(\alpha)=1$ for all $\alpha \in H'$ and consequently, 
$$(p^2-p)\sum_{\alpha \in H'}\overline{-(\phi(\alpha))}= -(p^2-p)|H'| = -(p^2-p)\frac{p^2-1}{d'}=-\frac{|G|}{d'}$$
Therefore, for cuspidal representation $\rho = \Theta_{\phi},$
$c_{\rho}= \begin{cases}
    0  \quad \text{when} \quad \phi^{d'}\neq 1\\
    -\frac{1}{d'} \quad \text{when} \quad \phi^{d'}= 1.
\end{cases}$

\section{Proof of the main theorem}\label{sec: proof of the main theorem}

\subsection{Zero Matrix}\label{subsec: zero matrix}
First, we figure out the case when $M=0$. 
\begin{align*}
    S_k(0) = \sum _{X \in \sH} \psi_0(X) =  \sum_{X\in \sH}\exp \left( \frac{2 \pi i }{p} \Tr(0) \right)= |\sH|.
\end{align*}
Therefore, for $M=0$ the contribution to \Cref{eq: N(C) formula} is 
\begin{equation}\label{eq: Main term M=0}
      \frac{1}{p^4} \psi_M(-C) S_k(0) S_k(0) = \frac{|\sH|^2}{p^4} =\left( \frac{1}{2}+\frac{1}{2d'}\right)^2\cdot p^4 + \sO(p^3)
\end{equation}

\subsection{Nonsingular Matrices}\label{subsec: nonsingular}
Recall  \cite[Section 1.2]{ganguly2023singular} For $M \in \GL(2, \FF_p)$ and $\rho \in \irr(\GL(2, \FF_p))$, 
\begin{equation}
G(\rho, M)= \rho(M)^{-1}G(\rho, Id)= g(\rho)\rho(M)^{-1}, g(\rho) ~\text{is a constant}.
\end{equation}
Hence $\Tr(G(\rho,M))= g(\rho)\chi_\rho(M)^{-1}.$ At this point we recall \cite[Theorem 1.5]{ganguly2023singular}.
\begin{theorem}\label{thm: Kondo's theorem} Suppose $\rho \in \irr(\GL(2, \FF_p))$, then
$$|g(\rho)|= \begin{cases}
    p \quad \text{when}~ \rho \simeq \mathrm{St} ~or~ \text{trivial} \\
    p^{3/2}\quad \text{when}~\rho \simeq I_{{\chi},1}, ~\chi ~\text{is nontrivial}\\
    p^2\quad  \text{otherwise}.
\end{cases}$$
\end{theorem}
Thus,  
$$|S_k(M)| \leq \sum_{\rho\in \irr(\GL(2, \FF_p))}|c_{\rho}||\Tr(G(\rho,M))| \leq \sum_{\rho\in \irr(\GL(2, \FF_p))}|c_{\rho}| |g(\rho)||\chi_\rho(M^{-1})| $$
\textbf{Case 1:} In this case, we assume that $B$ is scalar. Thus, $M$ and $MB$ belong to the same conjugacy class, so $|S_k(M)|= |S_k(MB)|$.
Hence,
\begin{align*}
    &|S_k(M)||S_k(MB)|= |S_k(M)||S_k(M)|=|S_k(M)|\overline{|S_k(M)|}\\
    &\leq  \sum_{\rho, \sigma \in \mathrm{Irr}(\GL(2, \FF_p))}|c_{\rho}|\overline{|c_\sigma|}|g_{\rho}|\overline{|g_\sigma|}|\chi_\rho(M^{-1})|\overline{|\chi_\sigma(M^{-1})|}
\end{align*}
This inequality along with character orthogonality yields
\begin{align*}
    \sum_{M \in \GL(2, \FF_p)}|S_k(M)||S_k(MB)|&\overset{\substack{\text{by character} \\ \text{orthogonality}}}{\leq} \sum_{\rho \in \mathrm{Irr}(\GL(2, \FF_p))}|c_{\rho}|^2|g_{\rho}|^2\sum_{M \in \GL(2, \FF_p)}|\chi_\rho(M^{-1})|\overline{|\chi_\rho(M^{-1})|} \\
    &= |G|\left(p^2|C_{1_G}|^2 + p^2|C_{St}|^2+ \sum_{\rho \in \mathrm{Cuspidal}}p^4|c_{\rho}|^2\right)\\
    &=|G|\left(p^2\left(\frac{1}{2}+\frac{1}{2d'}\right)+p^2\left(\frac{1}{2}-\frac{1}{2d'}\right)+p^4\left(\frac{d'-1}{2d'^2}\right)\right)\\
    &= \left(\frac{d'-1}{2d'^2}\right)\cdot p^8 + \sO(p^6)
\end{align*}
Thus, when $B$ is scalar,
\begin{equation}\label{eq: nonsingular count for scalar B}
     \sum_{M \in \GL(2, \FF_p)}|S_k(M)||S_k(MB)|= \left(\frac{d'-1}{2d'^2}\right)\cdot p^8 + \sO(p^6).
\end{equation}
\textbf{Case 2:} In this case, $B$ is non-scalar.
  \begin{align*}
     |S_k(M)| &\leq \sum_{\rho \in 1\mathrm{-dim}}|c_{\rho}||g(\rho)||\chi_\rho(M^{-1})| +  \sum_{\rho \in \mathrm{St \otimes \alpha}}|c_{\rho}||g(\rho)||\chi_\rho(M^{-1})| \\
     &+\sum_{\substack{\rho \in  \mathrm{Principal}\\ \mathrm{Series}}}|c_{\rho}||g(\rho)||\chi_\rho(M^{-1})|+\sum_{\rho \in \mathrm{Cuspidal}}|c_{\rho}||g(\rho)||\chi_\rho(M^{-1})|\\&\leq p\left(\frac{1}{2}+\frac{1}{2d'}\right)+p\left(\frac{1}{2}-\frac{1}{2d'}\right)|\chi_{St}(M^{-1})|+ {\frac{d' -1}{2d'}}\cdot p^2|\chi_{cusp}(M^{-1})|
\end{align*}
Note that $M$ and $MB$ cannot both be scalar simultaneously, by our choice of $B$. Therefore, 
\begin{align*}
    &|S_k(M)||S_k(MB)| \leq  \left({\frac{d' -1}{2d'}}\right)^2\cdot p^4|\chi_{cusp}(M^{-1})||\chi_{cusp}((MB)^{-1})|\\  
    &+\left({\frac{d' -1}{2d'}}\right)^2\cdot p^3\left(|\chi_{St}(M^{-1})||\chi_{cusp}((MB)^{-1})| + |\chi_{cusp}(M^{-1})||\chi_{St}((MB)^{-1})|\right) \\
    &+\left({\frac{d'^2 -1}{4d'^2}}\right)\cdot p^3\left(|\chi_{cusp}(M^{-1})| +|\chi_{cusp}((MB)^{-1})|\right) +\left({\frac{d' -1}{2d'}}\right)^2\cdot p^2|\chi_{St}(M^{-1})||\chi_{St}((MB)^{-1})|\\
    &+\left({\frac{d'^2 -1}{4d'^2}}\right)\cdot p^2\left(|\chi_{St}(M^{-1})| +|\chi_{St}((MB)^{-1})|\right) +\left({\frac{d' +1}{2d'}}\right)^2\cdot p^2 
\end{align*}
We use the character values of cuspidal representations, \Cref{lem: elliptic M goes to scalar MB iff B is nonsplit ss}, \Cref{lem: nonsplit ss M goes to nonsplit ss MB} and \Cref{lem: nonsplit ss M goes to non-semisimple MB} to simplify the first term and sum over all non-singular $M.$
\begin{align*}
    &\sA_1:=\left({\frac{d' -1}{2d'}}\right)^2\cdot p^4 \sum_{M \in \GL(2, \FF_p)}|\chi_{cusp}(M^{-1})||\chi_{cusp}((MB)^{-1})|\\
    &\leq \left({\frac{d' -1}{2d'}}\right)^2\cdot p^4 (\underbrace{\overbrace{2(p-1)^2}^{MB ~\text{non-scalar}}}_{M ~\text{scalar}}+ \underbrace{\overbrace{(p-1)^2}^{MB ~~\text{scalar}}+ \overbrace{2(p-1)(p^2-2)}^{MB ~\text{non-scalar}} }_{M ~\text{non-semisimple}}\\ &+ \underbrace{\overbrace{2(p-1)^2}^{MB ~~\text{scalar}}+\overbrace{(p+1)(p-1)^2}^{MB ~~\text{non-semisimple}} +\overbrace{p^4-1}^{MB ~~\text{elliptic}} } _{M ~\text{elliptic}} )\\
    &= \left({\frac{d' -1}{2d'}}\right)^2\cdot p^8 + \sO(p^7)
\end{align*}
Next, we simplify the second term by the character values of Steinberg and cuspidal representations and the three choices based on whether a scalar or non-scalar $M$ produces a scalar or non-scalar $MB$. 
\begin{align*}
  \sA_2&:= \left(\frac{d' -1}{2d'}\right)^2\cdot p^3\sum_{M \in \GL(2, \FF_p)}\left(|\chi_{St}(M^{-1})||\chi_{cusp}((MB)^{-1})| + |\chi_{cusp}(M^{-1})||\chi_{St}((MB)^{-1})|\right)\\
  &\leq \left(\frac{d' -1}{2d'}\right)^2\cdot p^3 \left(\underbrace{2p(p-1)}_{\substack{M~ \text{scalar}\\ MB ~\text{non-scalar}}}+\underbrace{(p-1)^2}_{\substack{M~ \text{non-scalar}\\ MB ~\text{scalar}}}+ \underbrace{2(p-1)(p^3-p-2)}_{\substack{M~ \text{non-scalar}\\ MB ~\text{non-scalar}}} \right)\\
   &+\left(\frac{d' -1}{2d'}\right)^2\cdot p^3 \left(\underbrace{(p-1)^2}_{\substack{M~ \text{scalar}\\ MB ~\text{non-scalar}}}+\underbrace{2p(p-1)}_{\substack{M~ \text{non-scalar}\\ MB ~\text{scalar}}}+ \underbrace{2(p-1)(p^3-p-2)}_{\substack{M~ \text{non-scalar}\\ MB ~\text{non-scalar}}} \right)\\
   &= \sO(p^7).
\end{align*} 
Next, the third term: 
\begin{align*}
    \sA_3&:= \left({\frac{d'^2 -1}{4d'^2}}\right)\cdot p^3\sum_{M \in \GL(2, \FF_p)}\left(|\chi_{cusp}(M^{-1})| +|\chi_{cusp}((MB)^{-1})|\right) \\
    &\leq \left({\frac{d'^2 -1}{4d'^2}}\right)\cdot p^3 \left(\underbrace{(p-1)(p+1)}_{\substack{M ~ \text{scalar}\\MB~\text{non-scalar}}} + \underbrace{(p-1)(p+1)}_{\substack{M~ \text{non-scalar}\\MB ~\text{scalar}}} + \underbrace{4(p-1)(p^3-p-2)}_{\substack{M~ \text{non-scalar}\\ MB ~\text{non-scalar}}} \right)\\
    &= \sO(p^7).
\end{align*}
Finally, the fourth term :
\begin{align*}
    \sA_4&:= \left({\frac{d' -1}{2d'}}\right)^2\cdot p^2 \sum_{M \in \GL(2, \FF_p)}|\chi_{St}(M^{-1})||\chi_{St}((MB)^{-1})|\\
    & \leq \left({\frac{d' -1}{2d'}}\right)^2\cdot p^2 \left(\underbrace{p(p-1)}_{\substack{M~ \text{scalar}\\ MB ~\text{non-scalar}}}+\underbrace{p(p-1)}_{\substack{M~ \text{non-scalar}\\ MB ~\text{scalar}}}+ \underbrace{(p-1)(p^3-p-2)}_{\substack{M~ \text{non-scalar}\\ MB ~\text{non-scalar}}} \right)\\
    &=\sO(p^6)
\end{align*}
The remaining terms are $\sO(p^6)$ and $\sO(p^2)$. Therefore, we have 
\begin{equation}\label{eq: nonsingular count}
    \sum_{M \in \GL(2, \FF_p)}|S_k(M)||S_k(MB)|\leq \left(\frac{1}{2}-\frac{1}{2d'}\right)^2\cdot p^8 + \sO(p^7)
\end{equation}

\subsection{Rank one Matrices}\label{subsec:rank one}
 For a rank one matrix $M \in \M(2,p)$, we  give a bound for $S_k(M)$.
 Using \cref{Eq:S_k(M)}, we have 
 \begin{align*}
     |S_k(M)| &\leq \sum_{\rho\in \irr(\GL(2, \FF_p))}|c_{\rho}||\Tr(G(\rho),M)| \\
     &= \sum_{\rho \in 1\mathrm{-dim}}|c_{\rho}||\Tr(G(\rho,M))| +  \sum_{\rho \in \mathrm{St \otimes \alpha}}|c_{\rho}||\Tr(G(\rho,M))| \\&+\sum_{\substack{\rho \in  \mathrm{Principal}\\ \mathrm{Series}}}|c_{\rho}||\Tr(G(\rho,M))| +\sum_{\rho \in \mathrm{Cuspidal}}|c_{\rho}||\Tr(G(\rho,M))|
\end{align*}
By \cite[Theorem 1.8]{ganguly2023singular}, $\Tr(G(\rho,M))=0$ when $\rho$ is not a trivial representation $1_G$, Steinberg representation $\mathrm{St}$, or the principal series representation $I_{\chi,1}$, where $\chi$ is a nontrivial character of $\FF_p^{\times}$. Moreover, in the same theorem, Ganguly and Rajan also computed the values of $\Tr(G(\rho,M))$ for the remaining ones. For the trivial representation $1_G$ and any rank one matrix $A \in \M(2, \FF_p)$
\begin{align*}
    &\Tr(G(1_G,A))= -p(p-1)
    \end{align*}
    Any rank one matrix in $\M(2,p)$ conjugate to either $A_a= \begin{pmatrix}
    a & 0 \\
    0 &0
\end{pmatrix}$, $a\neq 0$ or $N= \begin{pmatrix}
    0 &1\\
    0 &0
\end{pmatrix}$. For the Steinberg representation :
\begin{align*}
    &\Tr(G(\mathrm{St},A_a))= -p(p-1)\\
    &\Tr(G(\mathrm{St},N))= p^2(p-1)
\end{align*}
Finally, for the principal series representation $I_{\chi,1}$:
    \begin{align*}
    &\Tr(G(I_{\chi,1},A_a))= p(p-1)\overline{\chi(a)}G(\chi)\\
    &\Tr(G(I_{\chi,1},N))= p(p-1)G(\chi)
\end{align*}
where $G(\chi) $ is the classical Gauss sum.  
Substituting all the values of $c_\rho$ and $\Tr(G(\rho,M))$ into the $|S_k(M)|$ expression, we have the following two bounds for $M$ conjugate to $A_a$ and $N,$ respectively. 
\begin{align*}
   | S_k(M)| &< |C_{1_G}| |\Tr(G(1_G,A_a))| + |c_{\mathrm{St}}||\Tr(G(\mathrm{St},A_a))|\\
    &=p(p-1)\left(\frac{1}{2}+\frac{1}{2d'}\right)+p(p-1)\left(\frac{1}{2}-\frac{1}{2d'}\right)\\
    & = p(p-1)
\end{align*}
\begin{align*}
     | S_k(M)| &< |c_{1_G}| |\Tr(G(1_G,N))| + |c_{\mathrm{St}}||\Tr(G(\mathrm{St},N))|\\
    &=p(p-1)\left(\frac{1}{2}+\frac{1}{2d'}\right)+p^2(p-1)\left(\frac{1}{2}-\frac{1}{2d'}\right)\\
\end{align*}
\textbf{Case 1:} In this case, we assume that $B$ is scalar. Note that $B$ and $MB$ belong to same conjugacy class.\\ 
\begin{align*}
    &\sum_{0 \neq M \in \M(2,\FF_p)-\GL(2, \FF_p)}|S_k(M)||S_k(MB)|=\sum_{ M \sim A_a}|S_k(M)||S_k(MB)|+\sum_{ M \sim N}|S_k(M)||S_k(MB)|\\
    &\leq p^3(p-1)^3(p+1)+\left(p(p-1)\left(\frac{1}{2}+\frac{1}{2d'}\right)+p^2(p-1)\left(\frac{1}{2}-\frac{1}{2d'}\right)\right)^2(p^2-1)\\
    &=\left(\frac{1}{2}-\frac{1}{2d'}\right)^2\cdot p^8 + \sO(p^7)
\end{align*} 
Thus, for $B$ scalar,
\begin{equation}\label{eq: singular count for B scalar}
    \sum_{0 \neq M \in \M(2,\FF_p)-\GL(2, \FF_p)}|S_k(M)||S_k(MB)|=\left(\frac{1}{2}-\frac{1}{2d'}\right)^2\cdot p^8 + \sO(p^7)
\end{equation}
\textbf{Case 2:} In this case, $B$ is non-scalar.
\begin{align*}
    &\sum_{0 \neq M \in \M(2,\FF_p)-\GL(2, \FF_p)}|S_k(M)||S_k(MB)|=\sum_{ M \sim A_a}|S_k(M)||S_k(MB)|+\sum_{ M \sim N}|S_k(M)||S_k(MB)|\\
    &\leq \sum_{ \substack{M \sim A_a \\ MB \sim N}}|S_k(M)||S_k(MB)|+\sum_{ \substack{M \sim A_a \\ MB \sim A_a} }|S_k(M)||S_k(MB)|\\
    &+\sum_{  \substack{M \sim N\\ MB \sim A_a}}|S_k(M)||S_k(MB)|+\sum_{ \substack{M \sim N \\ MB \sim N} }|S_k(M)||S_k(MB)|\\
    &= (p^2-1)p(p-1) \left(p(p-1)\frac{d'+1}{2d'} + p^2(p-1)\frac{d'-1}{2d'}\right) +  (p^2+1)p^2(p-1)^3\\
    &+(p^2-1) p(p-1) \left(p(p-1)\frac{d'+1}{2d'} + p^2(p-1)\frac{d'-1}{2d'}\right) \\&+  2(p-1) \left(p(p-1)\frac{d'+1}{2d'} + p^2(p-1)\frac{d'-1}{2d'}\right)^2\\
    &= \sO(p^7)
\end{align*}
\begin{proof}[Proof of \Cref{thm: main result} ]
    For $B$ scalar \eqref{eq: nonsingular count for scalar B} and \eqref{eq: singular count for B scalar} obtain 
    \begin{align*}
       \frac{1}{p^4} \sum_{0 \neq M \in \M(2, \FF_p)}|S_k(M)||S_k(MB)|&\leq\left(\left(\frac{1}{2}-\frac{1}{2d'}\right)^2 +\frac{d'-1}{2d'^2}\right)\cdot p^4 + \sO(p^3)\\
        &= \frac{d'^2-1}{4d'^2} p^4 + \sO(p^3)
    \end{align*}
    The $M=0$ case \eqref{eq: Main term M=0} clearly dominates the last expression, guaranteeing a solution for $B$ scalar for sufficiently large $p$.
    Similarly, for $B$ non-scalar, \eqref{eq: nonsingular count} and Case $2$ of subsection \ref{subsec:rank one} yield
    \begin{align*}
        \frac{1}{p^4}\sum_{0 \neq M \in \M(2, \FF_p)}|S_k(M)||S_k(MB)|\leq \left(\frac{1}{2}-\frac{1}{2d'}\right)^2\cdot p^4 + \sO(p^3)
    \end{align*}
    which is also dominated by \eqref{eq: Main term M=0}, proving that $\sN(C)$ is positive for sufficiently large $p$.
    \end{proof}
\bigskip

\bibliographystyle{alpha}
\bibliography{References}
 \end{document}